\documentclass[]{spie}

\usepackage[]{graphicx}

\title{\Large \bf Quantum Algorithms for  the Jones Polynomial}

\author{Louis H. Kauffman\supit{a} and Samuel J. Lomonaco Jr.\supit{b}
\skiplinehalf
\supit{a} Department of Mathematics, Statistics and Computer Science  
(m/c 249), 851 South Morgan Street, University of Illinois at Chicago,
Chicago, Illinois 60607-7045, USA \\
\supit{b} Department of Computer Science and Electrical Engineering, University of
Maryland Baltimore County, 1000 Hilltop Circle, Baltimore, MD 21250, USA}
 
\authorinfo{Further author information: L.H.K.  E-mail: kauffman@uic.edu,
S.J.L. Jr.: E-mail: lomonaco@umbc.edu}
  
\begin{document} 

\newcommand{\Across}{\raisebox{-0.25\height}{\includegraphics[width=0.5cm]{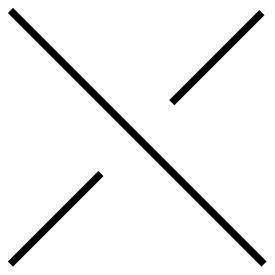}}}
\newcommand{\Bcross}{\raisebox{-0.25\height}{\includegraphics[width=0.5cm]{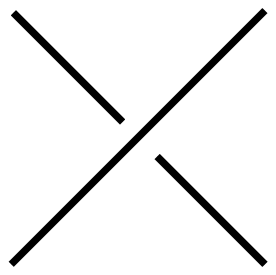}}}
\newcommand{\Asmooth}{\raisebox{-0.25\height}{\includegraphics[width=0.5cm]{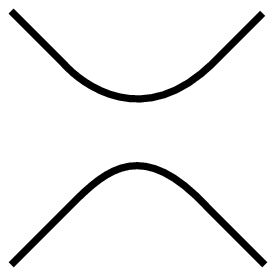}}}
\newcommand{\Bsmooth}{\raisebox{-0.25\height}{\includegraphics[width=0.5cm]{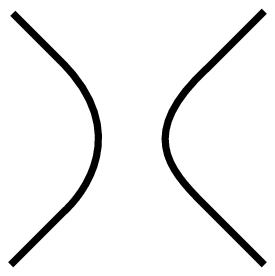}}}
\newcommand{\Rcurl}{\raisebox{-0.25\height}{\includegraphics[width=0.5cm]{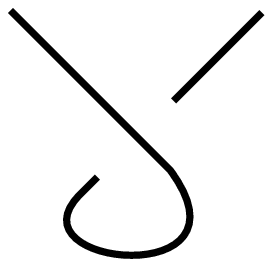}}}
\newcommand{\Lcurl}{\raisebox{-0.25\height}{\includegraphics[width=0.5cm]{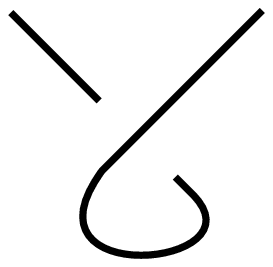}}}
\newcommand{\Arc}{\raisebox{-0.25\height}{\includegraphics[width=0.5cm]{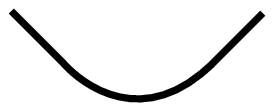}}}

 \maketitle

\begin{abstract}
This paper gives a generalization of the AJL algorithm for quantum computation of the Jones
polynomial to continuous ranges of values on the unit circle for the Jones parameter. We show 
that the Kauffman-Lomonaco 3-strand algorithm for the Jones polynomial is a special case of 
this generalization of the AJL algorithm.
\end{abstract}


\section{Introduction}
In \cite{QCJP} and in \cite{Three} we gave a quantum algorithm for computing the Jones polynomial via a unitary representation of the three-strand Artin braid group to the Temperley-Lieb algebra.
In the bracket polynomial version of this representation (see Section 2 of the present paper)
the representations were unitary for certain continuous ranges of choice of the polynomial variable
$A$ on the unit circle in the complex plane. In this paper we show that these three-strand representations are a subset of unitary representations of the Artin braid group on arbitrary numbers of strands and corresponding continuous ranges of the variable $A$ on the unit circle. These more 
general representations are in fact generalizations of the AJL representations \cite{Ah1,Ah2}
that were originally defined at certain roots of unity in the unit circle.
\bigbreak

The paper is organized as follows. In Section 2 we review the bracket polynomial model for the Jones polynomial, and its relationship with representations of the Temperley-Leib algebra. In Section 3 we review the 3-strand representation. In Section 4
we detail diagrammatically the construction of the generalized representation and show how it is 
related to the 3-strand representation and to the AJL representation. In Section 5 we give a 
diagrammatic proof of the requisite trace formula that is needed to make this representation into a quantum algorithm for computing the Jones polynomial. Much remains to be explored in these
directions. The present paper was sparked by our work in \cite{Group1} on NMR quantum computing, and there will be a sequel to the present paper \cite{Group2} that relates the present work to 
NMR research.
\bigbreak
 
\section{Bracket and Temperley Lieb Algebra}

The bracket polynomial \cite{KaB} model for the Jones polynomial \cite{JO,JO1,JO2,Witten} is usually described by the expansion
\begin{equation}
\langle \Bcross \rangle=A \langle \Bsmooth \rangle + A^{-1}\langle
\Asmooth \rangle \label{kabr}
\end{equation}

and we have

\begin{equation}
\langle K \, \bigcirc \rangle=(-A^{2} -A^{-2}) \langle K \rangle \label{kabr}
\end{equation}

\begin{equation}
\langle \Rcurl \rangle=(-A^{3}) \langle \Arc \rangle \label{kabr}
\end{equation}

\begin{equation}
\langle \Lcurl \rangle=(-A^{-3}) \langle \Arc \rangle \label{kabr}
\end{equation}
\bigbreak

\begin{figure}
     \begin{center}
     \begin{tabular}{c}
     \includegraphics[width=6cm]{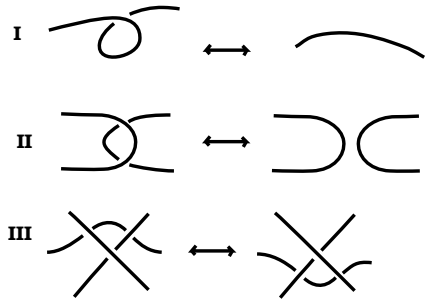}
     \end{tabular}
     \caption{\bf Reidemeister Moves}
     \label{Figure 1}
\end{center}
\end{figure}

The bracket expansion of a knot or link diagram is invariant under Reidemeister moves II and III
as shown in Figure 1, and can be normalized to be invariant under the first Reidemeister moves by 
multiplication by an appropriate power of $-A^{3}.$ Once normalized, it is a version of the Jones polynomial \cite{JO}, differing from it by a simple change of variable.    

The key idea behind the present quantum algorithms to compute the Jones polynomial is to use unitary 
representations of the braid group derived from Temperley-Lieb algebra representations that take the form
$$\rho(\sigma_{i}) = AI + A^{-1}U_{i}$$
where $\sigma_{i}$ is a standard generator of the Artin braid group, $A$ is a complex number of unit length, and
$U_{i}$ is a symmetric real matrix that is part of a representation of the Temperley-Lieb algebra. 
A diagrammatic version of the Temperley-Lieb algebra puts the form of this representation in exact
correspondence with the bracket expansion, where the parallel arcs 
$\Bsmooth$ correspond to the identity element of the algebra and the arcs in the form 
$\Asmooth$ correspond to the generator $U_{i}$ of the algebra when $\Bcross$ corresponds to the 
braid generator $\sigma_{i}.$ For more details about this
strategy and the background information about the Jones polynomial, the bracket model for the Jones polynomial and the Temperley-Lieb algebra
the reader may wish to consult \cite{Ah1,Ah2,JO,KaB,KA89,KL,KP,QCJP,Fibonacci,Spin,Three,QM,NSSF}. In the following sections, we have made use of such diagrammatic techniques and have included some material to make the paper partly self-contained.
\bigbreak

\section{Two Projectors and a Unitary Representation of the Three 
Strand Braid Group}
It is useful to think of the Temperley Lieb algebra as generated by projections
$e_{i} = U_{i}/\delta$ so that $e_{i}^{2} = e_{i}$ and $e_{i}e_{i\pm 
1}e_{i} = \tau e_{i}$ where
$\tau = \delta^{-2}$ and $e_{i}$ and $e_{j}$ commute for $|i-j|>1.$
\vspace{3mm}

With this in mind, consider elementary projectors
$e = |A \rangle \langle A|$ and $f=|B \rangle \langle B|$. We assume that $ \langle A|A \rangle = \langle B|B \rangle =1$ so 
that $e^{2} = e$ and $f^{2} =f.$
Now note that

$$efe = |A \rangle \langle A|B \rangle \langle B|A \rangle \langle A| = \langle A|B \rangle \langle B|A \rangle e = \tau e$$

\noindent Thus $$efe=\tau e$$

\noindent where $\tau = \langle A|B \rangle \langle B|A \rangle $.

This algebra of two projectors is the simplest instance of a 
representation of the Temperley Lieb
algebra.  In particular, this means that a representation of the 
three-strand braid group is
naturally associated with the algebra of two projectors.
\vspace{3mm}

Quite specifically if we let $\langle A| = (a,b)$ and $|A \rangle = (a,b)^{T}$ the 
transpose of this row vector,
then

$$e=|A \rangle \langle A| =  \left[
\begin{array}{cc}
      a^{2} & ab  \\
      ab & b^{2}
\end{array}
\right] $$

\noindent is a standard projector matrix when $a^{2} + b^{2} = 1.$ 
To obtain a specific
representation, let
\smallskip

l$e_{1} =  \left[
\begin{array}{cc}
      1 & 0  \\
      0 & 0
\end{array}
\right] $
and
$e_{2} =  \left[
\begin{array}{cc}
      a^{2} & ab  \\
      ab & b^{2}
\end{array}
\right] .$

\noindent  It is easy to check that
$e_{1}e_{2}e_{1} = a^{2}e_{1}$
and that
$e_{2}e_{1}e_{2} = a^{2}e_{2}.$

\noindent Note also that
$e_{1}e_{2}=  \left[
\begin{array}{cc}
      a^{2} & ab  \\
      0 & 0
\end{array}
\right]$
and
$e_{2}e_{1}=  \left[
\begin{array}{cc}
      a^{2} & 0  \\
      ab & 0
\end{array}
\right].$

We define $U_{i} = \delta e_{i}$ \, for $i=1,2$ with $a^{2}  = \delta^{-2}.$
Then we have , for $i = 1,2$
$$U_{i}^{2} = \delta U_{i} \, ,\, U_{1}U_{2}U_{1} = U_{1},\, U_{2}U_{1}U_{2} = U_{2}.$$
Thus we have a representation of the Temperley-Lieb algebra on three strands. See \cite{KL} for a discussion of the properties
of the Temperley-Lieb algebra.
\bigbreak

Note also that we have 
$$trace(U_{1})=trace(U_{2}) = \delta,$$
while
$$trace(U_{1}U_{2}) = trace(U_{2}U_{1}) = 1$$
where {\it trace} denotes the usual matrix trace.
These formulas show that the trace is working correctly with respect to the bracket evaluation.
See the last part of the present paper for a more extended discussion of this point.
\vspace{3mm}

\noindent Now we return to the matrix parameters:
Since $a^{2} + b^{2} = 1$ this means that $\delta^{-2} + b^{2} = 1$ whence
$b^{2} = 1-\delta^{-2}.$

\noindent Therefore $b$ is real when $\delta^{2}$ is greater than or 
equal to $1$.
\vspace{3mm}

We are interested in the case where $\delta = -A^{2} - A^{-2}$ and 
{\em $A$ is a unit complex
number}.  Under these circumstances the braid group representation

$$\rho(\sigma_{i}) = AI + A^{-1}U_{i}$$

\noindent will be unitary whenever $U_{i}$ is a real symmetric 
matrix. Thus we will obtain a
unitary representation of the three-strand braid group $B_{3}$ when 
$\delta^{2} \geq 1$.

\noindent For any $A$ with $d = -A^{2}-A^{-2}$ these formulas define a representation of the braid group. With 
$A=exp(i\theta)$, we have $d = -2cos(2\theta)$. We find a specific range of
angles $\theta$ in the following disjoint union of angular intervals
$$\theta \in [0,\pi/6]\sqcup[\pi/3,2\pi/3]\sqcup[5\pi/6,7\pi/6]\sqcup[4\pi/3,5\pi/3]\sqcup[11\pi/6,2\pi]$$
{\it that give unitary representations of
the three-strand braid group.} Thus a specialization of a more general represention of the braid group gives rise to a continuous family
of unitary representations of the braid group.
\bigbreak

\subsection{A Quantum Algorithm for the Jones Polynomial
on Three Strand Braids}

\label{threeStrands}

We gave above an example of a unitary representation of the three-strand 
braid group.In fact, we can use this representation to compute the Jones 
polynomial for closures of 3-braids, and therefore this
representation provides a test case for the corresponding quantum 
computation. We now analyse this case by first making
explicit how the bracket polynomial is computed from this 
representation. This unitary representation and its application to
a quantum algorithm first appeard in \cite{QCJP}. When coupled with 
the Hadamard test, this algorithm gets values for the Jones polynomial
in polynomial time in the same way as the AJL algorithm \cite{Ah1}. 
It remains to be seen how fast these algorithms are in principle when 
asked to compute the polynomial itself rather than certain specializations of it.
\vspace{3mm}

First recall that the representation depends on two matrices $U_{1}$ 
and $U_{2}$ with

$U_{1} =  \left[
\begin{array}{cc}
      \delta & 0  \\
           0 & 0
\end{array}
\right] $
and \,
$U_{2} =  \left[
\begin{array}{cc}
      \delta^{-1} & \sqrt{1-\delta^{-2}}  \\
          \sqrt{1-\delta^{-2}}  & \delta - \delta^{-1}
\end{array}
\right]. $

\noindent The representation is given on the two braid generators by

$$\rho(\sigma_{1})= AI + A^{-1}U_{1}$$
$$\rho(\sigma_{2})= AI + A^{-1}U_{2}$$

\noindent for any $A$ with $\delta = -A^{2}-A^{-2}$, and with
$A = exp(i\theta)$, then $\delta = -2cos(2\theta)$. We get the specific range of
angles $\theta \in [0,\pi/6]\sqcup[\pi/3,2\pi/3]\sqcup[5\pi/6,7\pi/6]\sqcup[4\pi/3,5\pi/3]\sqcup[11\pi/6,2\pi]$ that 
give unitary representations of
the three-strand braid group.
\vspace{3mm}

Note that $tr(U_{1})=tr(U_{2})= \delta$ while $tr(U_{1}U_{2}) = 
tr(U_{2}U_{1}) =1.$
If $b$ is any braid, let $I(b)$ denote the sum of the exponents in 
the braid word that expresses $b$.
For $b$ a three-strand braid, it follows that
$$\rho(b) = A^{I(b)}I + \tau(b)$$
\noindent where $I$ is the $ 2 \times 2$ identity matrix and 
$\tau(b)$ is a sum of products in the Temperley Lieb algebra
involving $U_{1}$ and $U_{2}.$ Since the Temperley Lieb algebra in 
this dimension is generated by $I$,$U_{1}$, $U_{2}$,
$U_{1}U_{2}$ and $U_{2}U_{1}$, it follows that
$$\langle \overline{b} \rangle = A^{I(b)}\delta^{2} + tr(\tau(b))$$
\noindent where $\overline{b}$ denotes the standard braid closure of 
$b$, and the sharp brackets denote the bracket polynomial
as described in previous sections. From this we see at once that
$$ \langle \overline{b} \rangle = tr(\rho(b)) + A^{I(b)}(\delta^{2} -2).$$
It follows from this calculation that the question of computing the 
bracket polynomial for the closure of the three-strand
braid $b$ is mathematically equivalent to the problem of computing 
the trace of the matrix $\rho(b).$
\bigbreak

The matrix in question is a product of unitary matrices, the quantum 
gates that we have associated with the braids
$\sigma_{1}$ and $\sigma_{2}.$ The entries of the matrix $\rho(b)$ 
are the results of preparation and detection for the
two dimensional basis of qubits for our machine:
$$ \langle i|\rho(b)|j \rangle.$$
\noindent Given that the computer is prepared in $|j \rangle $, the 
probability of observing it in state $|i \rangle $ is equal
to $| \langle i|\rho(b)|j \rangle|^{2}.$ Thus we can, by running the quantum 
computation repeatedly, estimate the absolute squares of the
entries of the matrix $\rho(b).$  This will not yield the complex 
phase information that is needed for either the trace of the
matrix or the absolute value of that trace.
\bigbreak
However, we do know how to write a quantum algorithm to compute the 
trace of a unitary matrix (via the Hadamard test).
Since $\rho(b)$ is unitary, we can use this approach to approximate 
the trace of $\rho(b).$ This yields a quantum algorithim for
the Jones polynomial for three-stand braids (evaluated at points $A$ 
such that the representation is unitary).
Knowing $tr(\rho(b))$ from the quantum computation, we then have the 
formula for the bracket, as above,
$$ \langle \overline{b} \rangle = trace(\rho(b)) + A^{I(b)}(\delta^{2} -2).$$
Then the normalized polynomial, invariant under all three 
Reidemeister moves is given by
$$f(\overline{b}) = (-A^{3})^{-I(b)} \langle \overline{b} \rangle .$$
Finally the Jones polynomial in its usual form is given by the formula
$$V(\overline{b})(t) = f(\overline{b})(t^{-1/4}).$$

Thus we conclude that our quantum computer can approximate values of 
the Jones polynomial.

\section{Generalizing the  AJL Representation}
\label{KL-AJL}
 
In this section we show how the KL (Kauffman-Lomonaco) algorithm described in the previous section
becomes a special case of a generalization of the AJL algorithm:
Here we use notation from the AJL paper. In that paper, the generators
$U_{i}$ (in our previous notation) for the Temperley-Lieb algebra, are denoted by
$E_{i}.$ We will first describe the AJL representation of the Temperley-Lieb Algebra and we will show how that representation works for a continuous range of values of the parameter $\theta$ described below. In the original treatment of AJK \cite{Ah1} the values are restricted to a discrete range corresponding to $exp(i \theta/2)$ a root of unity. We observed this phenomena of continuous ranges of values in our work on the three-strand model, described in the previous section, and found a way to extend it to AJL, as will be detailed below. This work has benefited from interaction with all of the authors of the paper \cite{Group1} and our collaboration in bringing these algorithms to application in NMR quantum computing. The contents of this section will be connected with NMR quantum computing in a paper that is under preparation \cite{Group2}.
\bigbreak

In this section we will construct matrix representations of the Temperley-Lieb algebra.
In the next section, we will discuss the structure of trace functions on these representations that can be used to produce quantum algorithms.
\bigbreak

Let $ \lambda_{k} = sin(k \theta).$ For the time being $\theta$ is an
arbitrary angle. Let $A = iexp(i \theta/2)$ so that
$d = -A^{2} - A^{-2} = 2cos(\theta).$
\bigbreak

We need to choose $\theta$ so that $sin(k \theta)$ is
non-negative for the range of $k$'s we use (these depend on the choice of
line graph as in AJL). And we insist that $sin(k \theta)$ is non-zero except
for $k=0.$ Then it follows from trigonometry that for all $k$
$$d = ( \lambda_{k-1} +  \lambda_{k+1})/ \lambda_{k}.$$ 
We shall see, below, that the values of $k$ will range from $0$ to 
$r$ for a fixed $r > 2$ in a given representation. We ask that $sin(r \theta)$ be greater than zero and
can take $\theta$ in the continuous range $0 < \theta \le \pi/r.$
\bigbreak

The AJL representation of the Temperley-Lieb algebra is based on the complex vector space
$H_{n}$ whose basis is  $\{ |i \rangle \}$ where $i$ is a ($0$,  $1$) bitstring of length $n.$ Each 
bitstring $i$ is seen as corresponding to a walk on a line-graph $G_{r}$ with $r-2$ edges and 
$r-1$ nodes. For example, the graph below is $G_{5}.$
$$1-----2-----3-----4$$
Bitstrings represent walks on a line graph with $0$ corresponding to a step to the left and $1$ corresponding to a step to the right. The walk begins at the left-most node, labelled 1. Thus $1011$ represents the walk
Right, Left, Right, Right ending at node number $3$ in the line graph above.
We say that a string of $n$ bits is a {\it walk on $G_{r}$} if  the walk remains inside $G_{r}$
for all its steps (i.e. one never faces the instruction to go beyond the right-most or left-most
end points of the graph). We let $H_{n,r}$ denote the subspace of $H_{n}$ spanned by bitstrings
that are walks on $G_{r}.$
\bigbreak

The representation of the Temperley-Lieb algebra is denoted by
$$\Phi: TL_{n} \longrightarrow Matr(H_{n,r})$$
where $TL_{n}$ denotes the $n$-strand Temperley-Lie algebra generated by
$I, E_{1},\cdots, E_{n-1}$ and $Matr(H_{n,r})$ denotes matrix mappings of this space in the 
given bitstring basis. 
The AJL representation is
given in terms of $E_{i}$ such that $E_{i}^2 = dE_{i}$ and the $E_{i}$ satisfy
the Temperley-Lieb relations. Each $E_{i}$ acts non-trivially at the $i$ and
$i+1$ places in the bit-string basis for the space. If $p$ denotes the bitstring and we are computing
$\Phi(|p \rangle),$ then we define  $ z(i)$ to be  the endpoint of the  walk
described by the bitstring $p$ using only the first $( i - 1)$ bits of $p$.  Each $E_{i}(p)$ is based
upon $ \lambda_{z(i)-1}, \lambda_{z(i)}, \lambda_{z(i)+1}.$  
In our example we have  $p = 1011$ represents the walk
Right, Left, Right, Right ending at node number $3$ in
$$1-----2-----3-----4,$$
and $z(1) = 1, z(2) = 2, z(3) = 1, z(4) = 1, z(5) = 3.$
\bigbreak

More precisely, if we let
$$|v(a) \rangle = (\sqrt{ \lambda_{a-1}/ \lambda_{a}}, \sqrt{ \lambda_{a+1}/ \lambda_{a}})^{T}$$
(i.e. this is a column vector. T denotes transpose.) Then
$$E_{i} = |v(z(i)) \rangle \langle v(z(i))|.$$
Here it is understood that this refers to the action
on the bitstrings $$----------01----------$$ and $$----------10----------$$ obtained from
the given bitstring by modifying the $i$ and $i+1$ places. The basis order is
$01$ before $10.$ \bigbreak

The explicit form of the the transformations $E_{i}$ is given by the equations below for a generic 
matrix $E.$ We need to explicate $\Phi(E_{i})$ and we shall simply write $E_{i}$ instead
of $\Phi(E_{i}).$

$$E =  v v^{T} = \left[
\begin{array}{cc}
      \lambda_{-}/\lambda_{0}  & \frac{ \sqrt{\lambda_{-} \lambda_{+}}}{\lambda_{0}}  \\
     \frac{ \sqrt{\lambda_{-} \lambda_{+}}}{\lambda_{0}}   &  \lambda_{+}/\lambda_{0}
\end{array}
\right] $$
Here $v^{T} = ( \sqrt{\frac{ \lambda_{-}}{\lambda_{0}}},  \sqrt{\frac{ \lambda_{+}}{\lambda_{0}}}).$
It is easy to see that $$E^{2} = (\frac{\lambda_{-} + \lambda_{+}}{\lambda_{0}}) E$$
since $E^{2} = v v^{T} v v ^{T} = (v^{T} v) E.$
For $E_{i}$, we take $\lambda_{-} = \lambda_{z(i) - 1},$
$\lambda_{+} = \lambda_{z(i) + 1}$ and $\lambda_{0} = \lambda_{z(i)}.$
For the action of $E_{i}(p)$ for a given bitstring $p,$ let $p|i$ denote the restriction of $p$ to the first 
$i-1$ bits in the string. Then we need to explicate $E_{i}$ at the bitstring $p$ and hence its values at 
$p|i\, k\, l\rangle$ where $k$ and $l$ are either $0$ or $1$. The transformation will not change bits beyond these two extra places. Then {\it by definition}
 $$E_{i}|p|i\, 0 \, 0 \rangle = E_{i}|p|i\, 1 \, 1 \rangle = 0$$ and $E_{i}$ acts on the remaining
 subspace according to the matrix formulas we have given above. This means that we can regard
 $E_{i}$ diagrammatically as a cup-cap combination. 
 $E = \raisebox{-0.40\height}{\includegraphics[width=.5 cm]{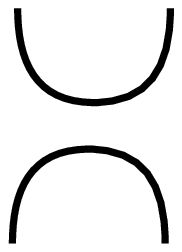}}$ where the legs of the diagram
 correspond to the indices that can be either $0$ or $1$. In this formalism 
$ \raisebox{-0.40\height}{\includegraphics[width=.5 cm]{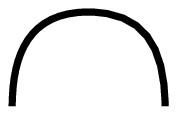}}$ takes the role of $v^{T}$ and 
$ \raisebox{-0.40\height}{\includegraphics[width=.5 cm]{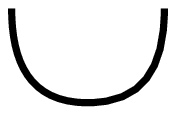}}$ the role of $v$ in the
decomposition $E = v v^{T}.$ Each cap or cup can receive only two indices and these must be different, since the transformations are zero when there is a repetition of $0$ or a repetition of $1.$
The only further relation that is needed to prove from this diagrammatic point of view is that 
$ \raisebox{-0.40\height}{\includegraphics[width=1.0 cm]{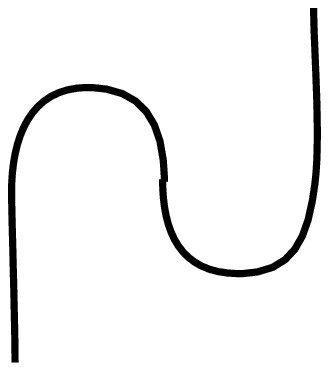}}$ is an identity transformation
as a mapping defined on a single bit. The relations $E_{i}E_{i+1}E_{i} = E_{i}$ follow from this.
To see that this composition is the identity, consider one of its cases as shown here:
$ \raisebox{-0.40\height}{\includegraphics[width=1.7 cm]{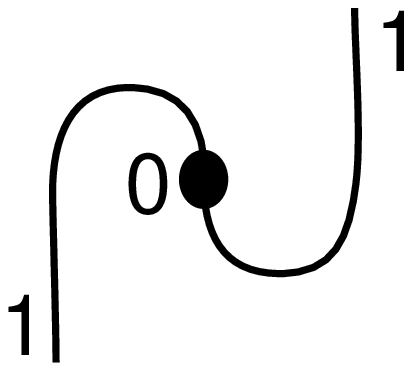}}.$ The fact that each cup
and each cap can only support a zero and a one, shows that the composition will necessarily be
a multiple of the identity. To see the details, we note that the local binary bit sequence must match
the indices, which in this case is $101.$ This means that $z(i+1) = z(i) + 1$ (due to the starting $1$ in 
the local bitstring, shifting to the next node in the graph).  The scalar contribution of 
 $ \raisebox{-0.40\height}{\includegraphics[width=1.7 cm]{couple01.eps}}$  is equal to 
 $$ S = \sqrt{  \frac{\lambda_{z(i)+1}}{\lambda_{z(i)}}    \frac{\lambda_{z(i+1)-1}}{\lambda_{z(i+1)}}    }$$
 where the first factor under the square root comes from the cap, and the second factor under the
 square root comes from the cup. The zero in the string $101$ corresponds to the $z(i+1)-1$
 in the second factor. Since we know that $z(i+1) = z(i) + 1$ it follows that $S=1$. This is one of the small
 number of cases to check that proves that 
 $ \raisebox{-0.40\height}{\includegraphics[width=1.0 cm]{couple.eps}}$ is the identity transformation.
 This completes the proof that $\Phi: TL_{n} \longrightarrow Matr(H_{n,r})$ is a representation of 
 the Temperley-Lieb algebra.
\bigbreak

\noindent {\bf Remark.} As we have seen in the discussion above, it helps  to take a diagrammatic point of view. Here we remark via Figure 2 on a simpler representation of the Temperley-Lieb algebra that is analogous to the AJL representation that we have discussed here. In Figure 2 we illustrate diagrams for the basic elements of the Temperley-Lieb algebra as we have explained them in the previous section. In order to build the 
elements $U_{i} = E_{i}$ we make them diagrammatically --  each a combination of a cup and a cap
(and appropriate identity lines), as illustrated in this figure. In the representation illustrated in this figure, each cup and each cap is represented by the same matrix $M$ and the necessary conditions for that $M$ are that the sum of the squares of its entries should be equal to $d = 2cos(\theta)$ (as above), and that the matrix product of the cup and cap matrices should be equal to the identity matrix. This is illustrated in the middle of Figure 2 by the identity showing the wavy line being pulled straight to form an identity line. All the diagrams in the figure correspond to matrices, with the indices for matrix elements corresponding to labelled endpoints of  the lines in the diagram. When two diagrammatic matrices are composed, an output line of one matrix is  attached to an input line of the other matrix. Thus in these diagrams, one sums over all the possible labels of an edge that has no free ends. This corresponds directly to the formula for matrix multiplication
where one sums the products of individual matrices over all double occurrences of indices.
\bigbreak

\begin{figure}
     \begin{center}
     \begin{tabular}{c}
     \includegraphics[width=6cm]{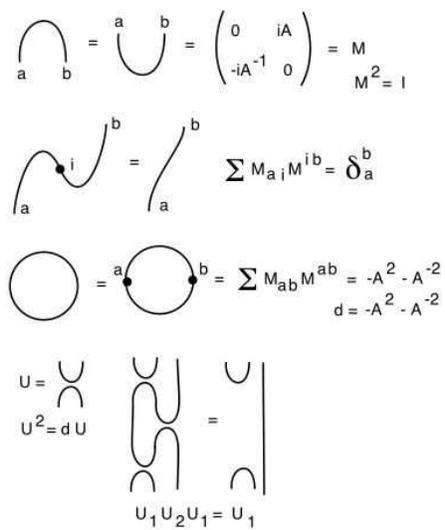}
     \end{tabular}
     \caption{\bf Diagrammatics for A Simple Temperley-Lieb Representation}
     \label{Figure 2}
\end{center}
\end{figure}

\begin{figure}
     \begin{center}
     \begin{tabular}{c}
     \includegraphics[width=5cm]{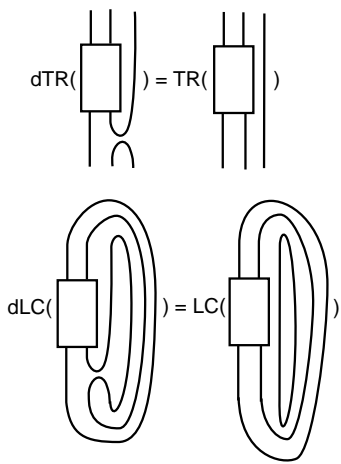}
     \end{tabular}
     \caption{\bf Trace Formula and Loop Count}
     \label{Figure 3}
\end{center}
\end{figure}

\begin{figure}
     \begin{center}
     \begin{tabular}{c}
     \includegraphics[width=10cm]{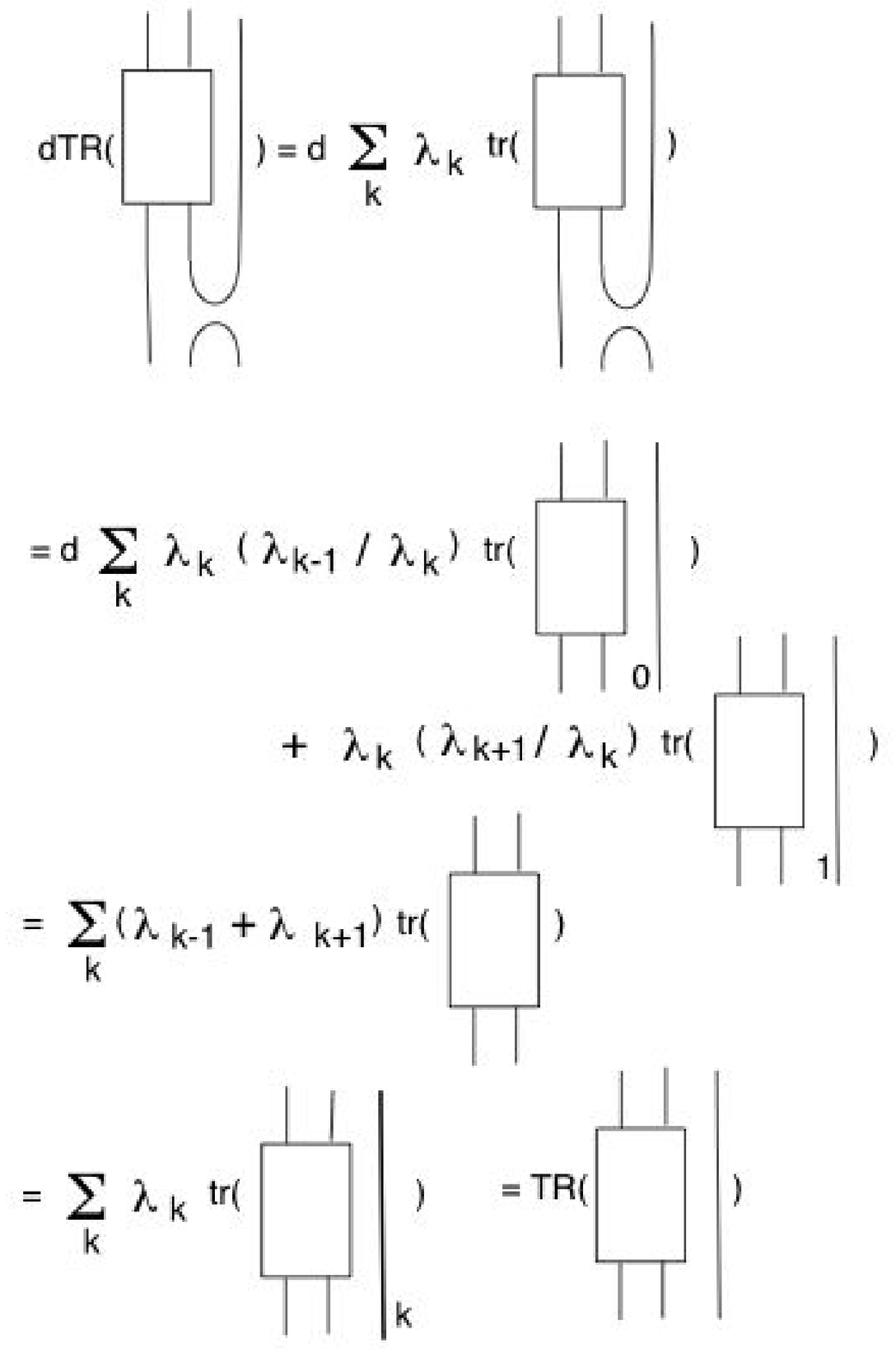}
     \end{tabular}
     \caption{\bf Proof of the Trace Formula}
     \label{Figure 4}
\end{center}
\end{figure}

\noindent {\bf Remark on the Three Strand Representation.} Now
look at the special case of a line graph with three nodes and two edges:
$$1-----2-----3.$$
The only admissible binary sequences are $|110 \rangle $ and $|101 \rangle ,$ so the space
corresponding to this graph is two dimensional, and it is acted on by
$E_{1}$ with $z(1) = 1$ in both cases
(the empty walk  terminates in the first node)
and
$E_{2}$ with $z(2) = 2$ for $|110 \rangle $ and $z(2)=2$ for $|101 \rangle .$
Then we have
$$E_{1}|110 \rangle = 0, E_{1}|101 \rangle = d |101\rangle ,$$
$$E_{2}|xyz \rangle = |v \rangle \langle v|xyz \rangle $$   (xyz = 101 or 110)
where $v = (\sqrt{1/d}, \sqrt{d - 1/d})^{T}.$
\bigbreak

If one compares this two dimensional representation of the three
strand Temperley - Lieb algebra and the corresponding braid group representation, with
the representation Kauffman and Lomonaco use in their paper, it is clear that it is
the same (up to the convenient replacement of $A = exp(i \theta)$ by $A = i exp(i \theta/2)$). The trace formula of AJL in this case is a variation of
the trace formula that Kauffman and Lomonaco use. See the next section for a general discussion of the trace. Note that the AJL algorithm as formulated in \cite{Ah1} does not
use the continuous range of angles that are available to the KL algorithm, but our generalization does 
allow this continuous angular range.
\bigbreak

\section {The Trace Function}
We now treat the trace function on the AJL representation from a diagrammatic point of view.
Let $\Phi: TL_{n} \longrightarrow Matr(H_{n,r})$ denote the representation of the Temperley-
Lieb algebra discussed in the last section.  Let $M = \Phi(\alpha)$ for any element $\alpha$ of the Temperley-Lieb algebra. We define a trace functional $TR(M)$ by the formula
$$Tr(M) = \sum_{k} \lambda_{k} tr(M_{k})$$ where $\lambda_{k} = sin(k \theta)$ as in the previous section,
$tr$ denotes standard matrix trace and $M_{k}$ denotes the restriction of $M$ to walks that end on the 
node $k$ in the graph $G_{r}.$ Just as in AJL, since this trace can be used to compute the bracket
polynomial of the closure of the braid at our admissible values of $A$ on the unit circle, it follows that there is a quantum algorithm for this computation by having the quantum computer evaluate the 
standard traces $tr(M_{k}).$ Each $M_{k}$ is a unitary matrix, and its trace can be found via the 
Hadamard test.
\bigbreak

We will prove that $$dTR(M E_{n}) = TR(M')$$ where $M'$ is the inclusion of $M$ in the corresponding
representation of $TL_{n+1}$ and $E_{n}$ denotes the matrix representation  of the element $E_{n}$ in $TL_{n+1}.$ This formula is exactly what is needed to have a Markov trace on the corresponding representation of the Artin braid group and hence to have a link invariant corresponding to the bracket model \cite{KaB} of the Jones polynomial. To see this from the diagrammatic point of view, examine
Figure 3. In that figure we have shown the diagrammatic version of the formula
$dTR(M E_{n}) = TR(M')$ at the top of the figure, and we have compared with the corresponding
loop count formula that corresponds to bracket polynomial calculation. The trace we define on 
the representation of the Temperley-Lieb algebra will agree with the loop counts performed by the 
bracket expansion on the closure of the braid, when the $TR$ formula is satisfied. In the figure
$LC$ refers to loop-count and denotes the evaluation of a collection of loops as $d^{c}$ where
$c$ is the number of loops in that collection.
\bigbreak

Finally, Figure 4 is a diagrammatic rendering of the proof of the trace formula. The expansion in the middle of the derivation corresponds to the evaluation of a representation of a Temperley-Lieb 
projector as described in the last section, and it takes into account the action of this projector on 
the bitstrings in the representation. In the middle of the calculation the labels $0$ and $1$ remind the reader of the bitstring values for which these terms are non-zero. At the end of the calculation, the label $k$ reminds the reader that these matrices are acting on walks that end on node $k$ in the graph. We have left some of the details to the reader, but the main line 
of the argument is in Figure 4. The reader should note that the syntax of this use of the diagrams is 
explained in the last section. This completes the proof of the trace formula and hence the proof that this generalization of the AJL algorithm to continuous angular ranges can be used as a quantum algorithm for the Jones polynonmial.
\bigbreak

\end{document}